\documentclass[10pt,twocolumn,twoside]{IEEEtran}

\usepackage{cite}

\ifCLASSINFOpdf
\else
\fi
\usepackage{array}

\usepackage{amssymb,amsmath}
\usepackage{graphicx}
\usepackage{amsthm}
\usepackage{epstopdf}
\allowdisplaybreaks[4]
\usepackage{amsmath}
\usepackage{flushend}
\usepackage{color}

\emergencystretch=\maxdimen
\newtheorem{theorem}{Theorem}
\newtheorem{definition}{Definition}

\newtheorem{lemma}{Lemma}

\newtheorem{assumption}{Assumption}
\newtheorem{remark}{Remark}

\newcommand{\col}{$\upshape{col}$}
\newcommand{\blk}{$\upshape{blk}$}
\newcommand{\diag}{$\upshape{diag}$}

\begin{document}
%
\title{Nash Equilibrium Seeking for General Linear Systems with Disturbance Rejection}
%
%
%

\author{Xin~Cai,
        Feng~Xiao,~\IEEEmembership{Member,~IEEE,}
        Bo~Wei,
        Mei~Yu,
        and Fang~Fang~\IEEEmembership{Senior Member,~IEEE,}
\thanks{This work has been submitted to the IEEE for possible publication. Copyright may be transferred without notice, after which this version may no longer be accessible. 
This work was supported by the National Natural Science Foundation of China (NSFC, Grant Nos. 61873074, 61903140). Corresponding author: F. Xiao.}

\thanks{X. Cai and F. Xiao are with the State Key Laboratory of Alternate Electrical Power System with Renewable Energy Sources and with the School of Control and Computer Engineering, North China Electric Power University, Beijing 102206, China (emails: caixin\_xd@126.com; fengxiao@ncepu.edu.cn). X. Cai is also with the School of Electrical Engineering, Xinjiang University, Urumqi 830047, China.}

\thanks{B. Wei, M. Yu and F. Fang are with the School of Control and Computer Engineering,
North China Electric Power University,
Beijing 102206, China (emails: bowei@ncepu.edu.cn; meiyu@ncepu.edu.cn; ffang@ncepu.edu.cn).}
}



%
%

\markboth{This work has been submitted to the IEEE for possible publication. Copyright may be transferred without notice}{}
%



\maketitle


\begin{abstract}
This paper explores aggregative games in a network of general linear systems subject to external disturbances. To deal with external disturbances, distributed strategy-updating rules based on internal model are proposed for the case with perfect and imperfect information, respectively. Different from existing algorithms based on gradient dynamics, by introducing the integral of gradient of cost functions on the basis of passive theory, the rules are proposed to force the strategies of all players to evolve to Nash equilibrium regardless the effect of disturbances. The convergence of the two strategy-updating rules is analyzed via Lyapunov stability theory, passive theory and singular perturbation theory. Simulations are presented to verify the obtained results.
\end{abstract}

\begin{IEEEkeywords}
aggregative games, external disturbances, linear systems, Nash equilibrium seeking.
\end{IEEEkeywords}

%
\IEEEpeerreviewmaketitle

\section{Introduction}
In recent years, the game theory has been widely applied in engineering community, such as power grids \cite{Hobbs.2007,Nekouei.2015,Liu.2018,Gharesifard.2016}, mobile sensor networks \cite{Stankovic.2012,Zhu.2013}, and communication networks \cite{Gharesifard.2013,Romano.2018}. A critical problem in the game theory is to seek Nash equilibrium (NE) which is the solution of game theory corresponding to the desired target of multi-agent systems.

The purpose of distributed Nash equilibrium seeking is to design the strategy-updating rule for players to reach the NE of the game, at which each agent has the minimum cost.  According to available information for agents, noncooperative games can be classified into the case with full/perfect information and the case with partial/imperfect information. In the full/perfect information case, early research was conducted under the assumption that every agent can observe the strategies of others\cite{Li.1987,Shamma.2005,Yin.2011}. In the setups of the existing results, a coordinator is usually required to collect global information and send it to all agents in a semi-decentralized framework\cite{Grammatico.2017,Persis.2020}. However, it is impractical for agents in large-scale networks to know actions of other agents due to agents' limited observation and communication capabilities. Furthermore, the semi-decentralized setting may be unreliable under attacks on the coordinator. Thus, an increasing interest is aroused in dealing with the issue of incomplete (imperfect) information games. To overcome the difficult that only partial information is available for agents, consensus protocols have been developed in multi-agent systems and they provide powerful tools for designing distributed discrete- and continuous-time NE seeking algorithms. Here we mainly focus on the continuous-time NE seeking
algorithms. In view of game models, a number of results have been developed on $n$-person
noncooperative games, aggregative games, two-network zero-sum games, $n$-coalition noncooperative games and differential games respectively\cite{Gharesifard.2016,Ye.2017,Gharesifard.2013,Ye.2018b,Vamvoudakis.2011}. In view of the classifications of communication graphs, there are some works
on the fixed undirected graphs, time-varying graphs, directed
graphs and switching graphs\cite{Gharesifard.2016,Liang.2017,Ye.2018,Deng.2019b},  where communication schemes mainly include continuous-time communications, event-triggered and self-triggered based interactions\cite{Gharesifard.2016,Cai.2020,Cortes.2016}. Additionally, to deal with local box constraints and nonsmooth cost functions, projected subgradient dynamics were proposed to seek the NE\cite{Lou.2016}. With the consideration of coupled constraints in networked scenarios, based on the primal-dual theory, gradient dynamics combined with the coordination of Lagrange multipliers were designed to seek the NE which satisfying the coupled constraints\cite{Lu.2019}.

A common assumption of all aforementioned literature is that the agent dynamics are first-order systems which do not completely describe the dynamical systems in real world. In practical applications of game theory, such as the coverage in mobile sensor networks, economical dispatch in electricity markets and consensus in multi-agent systems, it is desirable that multiple physical agents reach the NE of games corresponding to individual optimal strategies. Hence, it is of practical interest to study the noncooperative game theory for multi-agent systems with complex dynamics. Recent studies have
focused on the players with complex dynamics in the NE seeking, such as Euler-Lagrange systems \cite{Deng.2019}, networks of second-order systems \cite{Persis.2019}, multi-integrator agents \cite{Romano.2020, Bianchi.2019}, and nonlinear dynamic systems \cite{Zhang.2019}. In \cite{Deng.2019}, a control input of each agent, which is composed of a known nonlinear term and the gradient of its cost function, was designed to steer agents to the NE of some aggregative game. In \cite{Persis.2019}, a dynamic pricing mechanism was proposed to force passive nonlinear second-order systems to move to the predefined NE of noncooperative games with linear quadratic utility/cost functions. In \cite{Romano.2020}, distributed seeking algorithms composed of the damping, gradients of cost functions and observations of disturbances were designed for multi-integrator agents with disturbances. To remove the requirement of global knowledge of communication graphs, integral adaptive weights were introduced in the gradient-based NE seeking algorithm\cite{Bianchi.2019}. In \cite{Zhang.2019}, the strategy of each agent was steered by a nonlinear function and a control input which was the gradient of the involved cost function. The existing results on the NE seeking in noncooperative games with integrator-type agents primarily were given in the forms of gradient-based algorithms. However, approaches based on gradient dynamics cannot be directly used in noncooperative games with general linear systems, due to that various linear relationships among states may result in nonzero eigenvalues in the systems. Moreover, it is impractical to use the states of systems in the strategies interacted among agents. For complex physical systems, the outputs are more easily detected than the states of systems. How to regulate the outputs of general linear systems to reach the NE of noncooperative games is still a problem to be solved.

Aggregative games, as a special class of noncooperative games with cost functions depending on player's own strategy and the aggregation of strategies of all players, have attracted an increasing interest (see \cite{Deng.2019,Deng.2018,Deng.2019b,Liang.2017,Ye.2017,Shi.2019} and the references therein). Some work therein focuses on aggregative games with complex dynamic systems including Euler-Lagrange systems\cite{Deng.2019}, nonlinear systems\cite{Zhang.2019}. Furthermore, the competition among distributed energy resources in electricity market can be modeled as aggregative games \cite{Deng.2018,Deng.2019}. The output power of each generation system, as the strategy, is driven by the turbine-generator dynamics, which is a high-order linear system. These observations motivate us to study the aggregative games with high-order dynamics. To the best of our knowledge, the NE seeking of aggregative games with general linear systems has not been reported.

Arising from environment or communication, disturbances and noises are ubiquitous in practical dynamical systems. Although various NE seeking algorithms have been developed, these algorithms may fail to seek NE in the presence of disturbances. Recently, some work has focused on the distributed NE seeking of multi-agent systems with disturbances which may be generated by exosystems or just assumed to be bounded\cite{Romano.2020,Deng.2018,Zhang.2019,Ye.2020}. In \cite{Romano.2020,Deng.2018}, dynamic compensators were designed based on the internal models for integrator-type agents with the known models of disturbances. Zhang \textit{et~al.} extended the above work to the case with disturbances generated by exosystems with unknown parameters\cite{Zhang.2019}. Ye explored a NE seeking algorithm based on observers for players affected by bounded disturbances\cite{Ye.2020}. In this paper, we are trying to design observers to observe a class of deterministic disturbances modeled by exosystems.

The objective of this paper is to investigate the aggregative games of general linear systems with disturbances and design strategy-updating rules (NE seeking algorithms) for the players to seek the NE of the game. To the best of our knowledge, this is the first result on aggregative games of multiple general linear systems with disturbances.The contributions of this paper are given as follows.

1)  We study an aggregative game of all agents with
general linear dynamics. Different from the setups in which the internal states of agents are assumed to be strategies interacted among
agents in \cite{Deng.2019,Persis.2019,Romano.2020,Bianchi.2019,Zhang.2019}, the outputs of general linear systems are considered as the strategies of agents in this paper.

2)  A novel strategy-updating rule is proposed for general
linear systems with disturbances. Due to that the systems may have no zero eigenvalues, the existing gradient-based NE seeking algorithms are not applicable to the studied linear systems. To overcome this difficulty, the integrals of cost functions' gradients over time are introduced in the strategy-updating rule to force agents' outputs to arrive at the NE of aggregative games.

3)  The passivity theory plays a critical role in designing and
analyzing the proposed strategy-updating rule. However, in the existing literature,  Lyapunov stability theory is the main employed tool to analyze the convergence of NE seeking algorithms.

The organization of this paper is as follows. In Section \uppercase\expandafter{\romannumeral2}, some related preliminaries are provided. In Section \uppercase\expandafter{\romannumeral3}, the considered problem is formulated. In Section \uppercase\expandafter{\romannumeral4}, two strategy-updating rules are designed and their convergence is analyzed. In Section \uppercase\expandafter{\romannumeral5}, simulation examples are presented. The conclusions and future topics are stated in Section \uppercase\expandafter{\romannumeral6}.

\section{Preliminaries}

Notations: $\it{R}$ and $\it{R}^n$ denote the real numbers set and the $n$-dimensional Euclidean space, respectively. Given a vector $x \in\it{R}^n$, $\|x\|$ is the Euclidean norm. $\otimes$ denotes the Kronecker product. $A^T$ and $\|A\|$ are the transpose and the spectral norm of matrix $A$, respectively. Denote $\col(x_1,\ldots,x_n)=[x_1^T,\ldots,x_n^T]^T$. Given matrices $A_1$, \ldots, $A_n$, $\blk\{A_1,\cdots,A_n\}$ denotes the block diagonal matrix with $A_i$ on the diagonal. $I_n$ is the $n\times n$ identity matrix. $\mathbf 0$ denotes a zero matrix with an appropriate dimension. $\bf{1}_n$ and $\bf{0}_n$ are $n$-dimensional column vectors consisting of all 1s and 0s, respectively.

\subsection{Graph Theory}

Consider a graph $\mathcal{G}=(\mathcal{I},\mathcal{E})$, where $\mathcal{I}=\{1,\ldots,N\}$ is set of nodes, and $\mathcal{E}$ is set of edges. $\mathcal{G}$ is called an undirected graph, if $(i,j)\in \mathcal{E}$  implies $(j,i)\in \mathcal{E}$, for $\forall i,j\in \mathcal{I}$. $N_i \subset \mathcal{I}$ is the set of nodes directly connected to node $i$, and it is called the neighbor set of node $i$. $A=[a_{ij}]\in \it{R}^{N \times N}$ is the adjacency matrix of graph $\mathcal{G}$, where $a_{ij}=1$, if $(j,i)\in \mathcal{E}$ and $a_{ij}=0$ otherwise. $D=\diag\{d_1,\cdots,d_N\}$ is degree matrix, where $d_i=\sum_{j\not=i}a_{ij}$. $L=D-A\in \it{R}^{N \times N}$ is the Laplacian matrix of graph $\mathcal{G}$. A path between two nodes is an ordered sequence of nodes forming consecutive edges. If there is a path between any two nodes, the graph is connected. For an undirected and connected graph, $L$ is positive semi-definite and $L \bf{1}_N=\bf{0}_N$. Besides a simple zero eigenvalue, all other eigenvalues of $L$ are positive.

\subsection{Input-to-state Stability}

Consider the cascade system
\begin{subequations}
\begin{align}
\dot{x}_1&=f_1(x_1,x_2), \label{11}\\
\dot{x}_2&=f_2(x_2),  \label{12}
\end{align}
\end{subequations}
where $f_1:\it{R}^{n_1}\times\it{R}^{n_2}\rightarrow\it{R}^{n_1}$ and $f_2: \it{R}^{n_2}\rightarrow\it{R}^{n_2}$ are locally Lipschitz. For the system~(\ref{11}), with any initial state $x_1(t_0)$ and any bounded input $x_2(t)$, it said to be input-to-state stable \cite[Definition 4.7]{Khalil.2002}, if there exist a class $\mathcal{KL}$ function $\beta$ and a class $\mathcal{K}$ function $\gamma$ such that the solution $x_1(t)$ exists for all $t\geq t_0$ and satisfies
\begin{equation*}
\|x_1(t)\|\leq \beta(\|x_1(t_0)\|,t-t_0)+\gamma (\sup_{t_0\leq\tau\leq t}\|x_2(\tau)\|).
\end{equation*}

\begin{theorem}\label{thm3}
(\cite[Theorem 4.19]{Khalil.2002}) Let $V:\it{R}^{n_1}\rightarrow \it{R}$ be a continuously differentiable function such that
\begin{equation*}
\begin{aligned}
\alpha_1(\|x_1\|)&\leq V(x_1)\leq \alpha_2(\|x_1\|),\\
\frac{\partial V}{\partial x_1}f_1(x_1,x_2)&\leq-W(x_1), \ \ \forall \|x_1\|\geq \rho(\|x_2\|)>0
\end{aligned}
\end{equation*}
$\forall x_1\in \it{R}^{n_1}, x_2\in \it{R}^{n_2}$, where $\alpha_1$, $\alpha_2$ are class $\mathcal{K}_\infty$ functions, $\rho$ is a class $\mathcal{K}$ function, and $W(x_1)$ is a continuous positive definite function on $\it {R}^{n_1}$. Then, the system~(\ref{11}) is input-to-state stable with $\gamma=\alpha_1^{-1}\circ\alpha_2\circ\rho$.
\end{theorem}

\begin{lemma}\label{lemma1}
(\cite[Lemma 4.7]{Khalil.2002}) If the system~(\ref{11}), with $x_2$ as input, is input-to-state stable and the origin of~(\ref{12}) is globally uniformly asymptotically stable, then the origin of the cascade system~(\ref{11}) and~(\ref{12}) is globally uniformly asymptotically stable.
\end{lemma}

\subsection{Passive linear system}

Consider the linear time-invariant system
\begin{subequations} \label{lti}
\begin{align}
\dot{x}&=Ax+Bu,\\
y&=Cx,
\end{align}
\end{subequations}
where $x\in \it{R}^n$, $u\in \it{R}$ and $y\in \it{R}$ are the state, control input and output, respectively. $(A,B)$ is controllable and $(A,C)$ is observable. The transfer function of system \eqref{lti} is denoted by $G(s)=C(sI-A)^{-1}B$. If all poles of $G(s)$ have nonpositive real parts, $G(s)$ is called positive real \cite[Definition 6.4]{Khalil.2002}. Some notions about the passivity of linear time-invariant system \eqref{lti} are give as follows.

\begin{definition} \label{passive}
(\cite[Definition 6.3]{Khalil.2002}) The system \eqref{lti} is said to be passive if there exists a continuous differentiable positive semidefinite function $V(x)$ such that
\begin{equation*}
\dot{V}=\frac{\partial V}{\partial x}(Ax+Bu)\leq u^Ty,\ \forall (x,u)\in \it{R}^n\times \it{R}
\end{equation*}
\end{definition}

\begin{lemma} \label{lem6.2}
(\cite[Lemma 6.2]{Khalil.2002}) $G(s)$ is positive real if and only if there exist matrices $P=P^T>0$, $Q$ such that
\begin{equation*}
\begin{aligned}
PA+A^TP&=-QQ^T,\\
PB&=C^T.
\end{aligned}
\end{equation*}
\end{lemma}

\begin{lemma} \label{lem6.4}
(\cite[Lemma 6.4]{Khalil.2002}) The linear time-invariant minimal realization \eqref{lti} with $G(s)$ is passive if $G(s)$ is positive real.
\end{lemma}

\subsection{Aggregative Games}
Consider a game $G=(\mathcal{I},\Omega,J)$ with $N$ players denoted by $\mathcal{I}=\{1,\ldots,N\}$. $\Omega=\Omega_1\times\cdots\times\Omega_N\subset\it{R}^{Nn}$ where $\Omega_i$ is the action set of agent $i\in \mathcal{I}$. $J$ $=$ $(J_1,\ldots,J_N)$, where $J_i(y_i,y_{-i})$ $:$ $\Omega_i\times\prod_{j\neq i}\Omega_j\rightarrow \it{R}$ is the cost function of player $i$, $y_i$ $\in$ $\Omega_i$ $\subset$ $\it{R}^{n}$ is its strategy, $y_{-i}=\col(y_1,\ldots,y_{i-1},y_{i+1},\ldots,x_N)$ is the strategy vector of players except player $i$.  Let $\boldsymbol{y}=(y_i,y_{-i})\in \Omega$ denote the strategy profile of all players, which also implies that $(y_i,y_{-i})=(y_1,\ldots,y_N)$. Every player controls his own strategy to minimize his cost function for given $y_{-i}$, i.e., $\min_{y_i\in \Omega_{i}} J_i(y_i,y_{-i})$.

Here we consider an aggregative game (refer to \cite{Daron.2013,Jensen.2010}) defined as follows.

\begin{definition}
The game $G=(\mathcal{I},\Omega,\tilde{J})$ is said to be an aggregative game with aggregator $\sigma:\Omega\rightarrow \it{R}^n$  and a function $\tilde{J_i}: \Omega_i\times\it{R}^n \rightarrow \it{R}$ such that for each player $i$
\begin{equation}\label{cf1}
J_i(y_i,y_{-i})=\tilde{J_i}(y_i,\sigma(\boldsymbol{y})), \  \text{for all}~ \boldsymbol{y}\in \Omega,
\end{equation}
\end{definition}
\noindent where $\sigma(\boldsymbol{y})=\sum_{i=1}^N y_i$ denotes the aggregation of all players' decisions.

\begin{definition}
(\cite[Definition 3.7]{Basar.1999}) Let $G=(\mathcal{I},\Omega,\tilde{J})$ be an aggregative game. Then $\boldsymbol{y}^*$ $=$ $(y_i^*,y_{-i}^*)$ is a pure strategy Nash equilibrium if for each player $i\in \mathcal{I}$
\begin{equation*}
y_i^*=\arg\min_{y_i\in \it \Omega_i}\tilde{J_i}(y_i,\sigma(y_i,y_{-i}^*)).
\end{equation*}
\end{definition}

At the Nash equilibrium point, no player will unilaterally change his strategy for less cost. The following three assumptions are widely used in \cite{Gadjov.2019,Deng.2019,Deng.2018,Shi.2019}.

\begin{assumption} \label{as2}
For all $i\in\mathcal{I}$, $\Omega_i=\it{R}^{n}$, the cost function $J_i(y_i,y_{-i})$ is continuously differentiable in its arguments and convex in $y_i$ for every fixed $y_{-i}$.
\end{assumption}

By Corollary 4.2 in \cite{Basar.1999} and \cite{Gadjov.2019}, under Assumption \ref{as2}, the game $G=(\mathcal{I},\Omega,\tilde{J})$ exists a pure NE $\boldsymbol{y}^*\in\Omega$ satisfying
\begin{equation}\label{nec1}
\nabla_{y_i}{J_i(y_i^*,y_{-i}^*)}=\mathbf{0}_{n},   \  \forall{i}\in\mathcal{I},
\end{equation}
where $\nabla_{y_i}{J_i(y_i,y_{-i})}$ is the gradient of player $i$'s cost function with respect to his strategy $y_i$. Denote $\phi{(\boldsymbol{y})}=\col(\nabla_{y_1}{J_1(y_1,y_{-1})},\ldots,\nabla_{y_N}{J_N(y_N,y_{-N})})$ which is called the pseudo-gradient mapping. Thus, from \eqref{nec1}, we can get
\begin{equation} \label{psegc}
\phi(\boldsymbol{y}^*)=\mathbf{0}_{Nn}.
\end{equation}

For the game with imperfect information, each player only communicate with neighbors to estimate the aggregator. The estimation is denoted by $\eta_i \in\it{R}^n$. To prepare for our development, define the map $F_i:\Omega_i\times\it{R}^{n}\rightarrow\it{R}^{n}$ as
\begin{equation} \label{pseg}
\begin{split}
F_i(y_i,\eta_i)=(\nabla_{y_i}{\tilde{J_i}(y_i,\sigma)}+\nabla_\sigma{\tilde{J_i}(y_i,\sigma)})\mid_{\sigma=\eta_i}, \forall{i\in\mathcal{I}},
\end{split}
\end{equation}

Furthermore, denote $F(\boldsymbol{y},\eta)=\col(F_1(y_1,\eta_1),\ldots,$ $F_N(y_N, \eta_N))$, where $\eta$ $=$ $\col(\eta_1,\ldots,\eta_N)$. When the estimation $\eta_i=\sigma$,  $\forall i\in{\mathcal{I}}$, from \eqref{cf1} and \eqref{pseg}, it is observed that $\phi(\boldsymbol{y})=F(\boldsymbol{y},\sigma)$.

\begin{assumption} \label{as3}
 $\phi :\Omega\rightarrow\it{R}^{Nn}$ is strongly monotonic, i.e., $(\boldsymbol{y-y{'}})^T(\phi(\boldsymbol{y})-\phi(\boldsymbol{y{'}}))\geq \mu\|\boldsymbol{y-y{'}}\|^2$, for $\forall \boldsymbol{y}$, $\boldsymbol{y{'}}\in\Omega$ and $\mu>0$. Moreover, $F(\boldsymbol{y},\eta)$ is  Lipschitz continuous in $\eta\in\it{R}^{Nn}$ for every $\boldsymbol{y}\in\Omega$, i.e., there exists a positive constant $\theta$ such that $\|{F(\boldsymbol{y},\eta)-F(\boldsymbol{y},\eta{'})}\| \leq \theta\|{\eta-\eta{'}}\|$,  $\forall \eta$, $\eta{'}\in\it{R}^{Nn}$.
\end{assumption}

Under Assumption \ref{as3}, the aggregative game has a unique (pure) NE \cite[Theorem 3]{Scutari.2014}.

\section{Problem Formulation}

We consider that each agent (player) in the network can be modeled by the following single-input single-output (SISO) linear system
\begin{equation} \label{ad1}
\left\{\begin{array}{l}
\dot{x}_i=A_ix_i+B_i(u_i+d_i),\\
y_i=C_ix_i,
\end{array}\right.
\end{equation}
where $x_i\in\it{R}^{n}$, $u_i\in\it{R}$, and $y_i\in\it{R}$ are the agent $i$'s state, control input and output respectively. $A_i$, $B_i$ and $C_i$ are constant matrices with appropriate dimensions. Note that we focus on the case of one-dimensional output for ease of exposition. However, the forthcoming results also hold for higher dimensions. Assume that $(A_i, B_i)$ is controllable and $(A_i, C_i)$ is observable. $d_i\in\it{R}$ is the local disturbance generated by the exosystem:
\begin{equation} \label{dis}
\left\{\begin{array}{l}
\dot{\nu}_i=S_i\nu_i,\\
d_i=U_i\nu_i,
\end{array}\right.
\end{equation}
where $\nu_i\in\it{R}^{q}$ is the internal state of the exosystem. The matrix pair $(S_i,U_i)$ is observable. Suppose that, for each $i=1$, \ldots, $N$, all the eigenvalues of the matrix $S_i$ are distinctively lying on the imaginary axis, which means that the disturbance is bounded. Here, we deal with the effects of a class of deterministic disturbances on the system \eqref{ad1}. Therefore, it is assumed that agent $i$ knows its $S_i$ and $U_i$, but is ignorant of the internal state $\nu_i$, for all $i\in\mathcal{I}$.

Since agent $i$ is not able to broadcast all of its states, the state $x_i$ is unsuitable to be used by all players in the game. Therefore, we consider the output $y_i$ as the strategy of agent $i$ and $y_{-i}=(y_1,\ldots,y_{i-1},y_{i+1},\ldots,y_N)$ as the output vector of other players. Let $\boldsymbol{y}=(y_1,\ldots,y_N)$. The purpose of each agent is
\begin{equation*}
\min_{y_i\in \it{R}} J_i(y_i,y_{-i})=\tilde{J_i}(y_i,\sigma(\boldsymbol{y})),
\end{equation*}
where $\sigma(\boldsymbol{y})=\sum_{i=1}^N y_i$. Meanwhile, disturbances appearing in agent dynamics can affect its strategy. The objective in this paper is to design strategy-updating rule for every player with dynamics \eqref{ad1} to achieve the NE of the aggregative game and to reject the disturbances simultaneously.

\begin{remark}
In this paper, the proposed method is applied to steer agents' strategies to the unique Nash equilibrium of the aggregative game, but not to the Pareto optimality. If the assumptions on cost functions are modified to ensure that the game has multiple Nash equilibria, some incentive or punishment mechanisms need to be designed to force agents to change their strategies to move from an undesired Nash equilibrium to the desired one that is the Pareto optimality. The related topic is out of the scope of this paper and will be addressed in our future work.
\end{remark}

\section{Main Results}

In this section, for the perfect information games, a strategy-updating rule is first developed for the players with dynamics \eqref{ad1}. Then, we consider the case with imperfect (incomplete) information and propose an improved strategy-updating rule. Moreover, the convergence of two strategy-updating rules is analyzed.

\subsection{Strategy-updating Rule with Perfect Information}
\begin{figure}   
\begin{center}
\includegraphics[width=8cm]{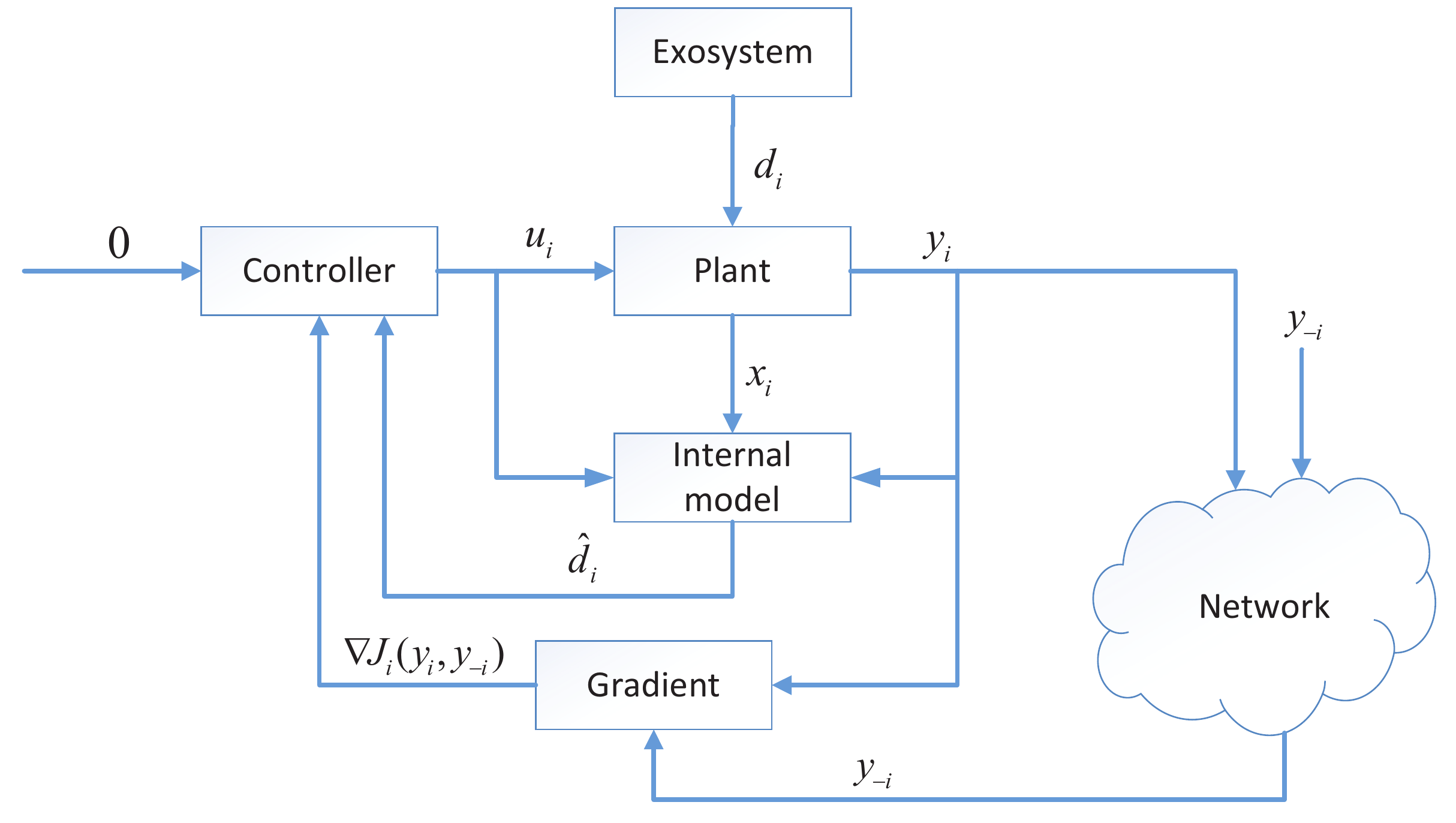}
\caption{Block diagram of the designed system with perfect information}
\label{fig1}
\end{center}
\end{figure}

In the game with perfect information, each player can obtain the strategies of other players, which can be considered as a special case of complete communication graph. For agents disturbed by external disturbances, combined with internal model, gradient and its integral, a strategy-updating rule  is proposed for agents to update strategies. It is shown in Fig. \ref{fig1}.

Let $e_i=-\nabla_i{J_i(y_i,y_{-i})}$ be the gradient of the player $i$'s cost function, $\gamma_i$ denote an auxiliary variable with $\dot{\gamma_i}=e_i$, and $\hat{d_i}$ be the observation of disturbances. The strategy-updating rule for player $i\in\mathcal{I}$ is designed as follows.
\begin{equation}\label{pi}
\begin{aligned}
\dot{x}_i&=(A_i-B_iK_i)x_i+B_i(k_{P_i}e_i+k_{I_i}\gamma_i-\hat{d_i}+d_i), \\
y_i &=C_ix_i,
\end{aligned}
\end{equation}
where $K_i\in\it{R}^{1\times n}$ such that $(A_i-B_iK_i)$ is Hurwitz,  $k_{P_i}$ and $k_{I_i}$ are positive constants.

To cope with external disturbances, the observer of disturbances is designed as follows based on the idea of internal model.
\begin{equation} \label{do}
\begin{aligned}
\hat{d}_i&=U_i(z_i+k_{o_i}B_i^Tx_i),\\
\dot{z}_i&=S_iz_i-k_{o_i}B_i^TB_i(k_{P_i}e_i+k_{I_i}\gamma_i)\\
         &\ \ \ +(S_ik_{o_i}B_i^T-k_{o_i}B_i^TA_i+k_{o_i}B_i^TB_iK_i)x_i,
\end{aligned}
\end{equation}
where $z_i$ is the internal state of the observer, and $k_{o_i}\in\it{R}^{q}$ satisfying that $(S_i-k_{o_i}B_i^TB_iU_i)$ is Hurwitz.

In order to analyze the properties of systems \eqref{pi} and \eqref{do}, let $\chi_i=\big[\begin{smallmatrix}x_i\\ \gamma_i\end{smallmatrix}\big]$ be the augmented state of the system,  $\rho_i=\nu_i-(z_i+k_{o_i}B_i^Tx_i)$ be the observation error, and $H_i=A_i-B_iK_i$. The closed-loop system \eqref{pi} and \eqref{do} can be rewritten as
\begin{equation}\label{red}
\begin{aligned}
\dot{\chi}_i&=\mathcal{A}_i\chi_i+\mathcal{B}_ie_i+\mathcal{D}_i\rho_i,\\
y_i&=\mathcal{C}_i\chi_i,\\
\dot{\rho}_i&=(S_i-k_{o_i}B_i^TB_iU_i)\rho_i,
\end{aligned}
\end{equation}
where
\begin{equation*}
\begin{aligned}
\mathcal{A}_i&=\begin{bmatrix}
H_i & B_ik_{Ii}\\
\boldsymbol{0}^T_n & 0
\end{bmatrix},\ \
\mathcal{B}_i=
\begin{bmatrix}
B_ik_{Pi}\\
1
\end{bmatrix},\\
\mathcal{C}_i&=
\begin{bmatrix}
C_i & 0
\end{bmatrix},\ \
\mathcal{D}_i=
\begin{bmatrix}
B_iU_i\\ \boldsymbol{0}_q^T
\end{bmatrix},
\end{aligned}
\end{equation*}
$e_i\in \it{R}$ is the control input, and $\rho_i\in \it{R}^q$ is the reference input. We will analysis the convergence of system \eqref{red} to show the effectiveness of the proposed method in the sequel.

First, we write system \eqref{red} in a stack form for convenience. Denote $\boldsymbol{\chi}=\col(\chi_1,\ldots,\chi_N)$, $\boldsymbol{e}=[e_1,\ldots,e_N]^T$, $\boldsymbol{\rho}=\col(\rho_1,\ldots,\rho_N)$ and $\boldsymbol{y}=[y_1,\ldots,y_N]^T$. Then, \eqref{red} can be expressed as
\begin{equation} \label{sf1}
\begin{aligned}
\boldsymbol{\dot{\chi}}&={\mathcal{A}}\boldsymbol{\chi}+\mathcal{B}\boldsymbol{e}+\mathcal{D}\boldsymbol{\rho},\\
\boldsymbol{y}&=\mathcal{C}\boldsymbol{\chi},\\
\boldsymbol{\dot{\rho}}&=(\boldsymbol{S}-K_o\boldsymbol{B}^T\boldsymbol{B}\boldsymbol{U})\boldsymbol{\rho},
\end{aligned}
\end{equation}
where $\mathcal{A}$ $=$ $\blk\{\mathcal{A}_1,\cdots,\mathcal{A}_N\}$, $\mathcal{B}$ $=$ $\blk\{\mathcal{B}_1,\cdots,\mathcal{B}_N\}$,
$\mathcal{C}$ $=$ $\blk\{\mathcal{C}_1,\cdots,\mathcal{C}_N\}$,
$\mathcal{D}$ $=$ $\blk\{\mathcal{D}_1,\cdots,\mathcal{D}_N\}$, $\boldsymbol{S}$ $=$ $\blk\{S_1,\cdots,S_N\}$,
$\boldsymbol{B}$ $=$ $\blk\{B_1,\cdots,B_N\}$, $\boldsymbol{U}$ $=$ $\blk\{U_1,\cdots,U_N\}$, and $K_o$ $=$ $\blk\{k_{o_1},\cdots,k_{o_N}\}$.

Next, the relationship between the equilibrium point of \eqref{sf1} and the NE of aggregative game $G=(\mathcal{I},\Omega,\tilde{J})$ is given as follows.

\begin{lemma} \label{lemma3}
Under Assumption \ref{as2}, if $(\boldsymbol{\chi^*}, \boldsymbol{0}_{Nq})$ is an equilibrium point of \eqref{sf1}, $\boldsymbol{y^*}$ is a NE of aggregative game $G=(\mathcal{I},\Omega,\tilde{J})$. Conversely, if $\boldsymbol{y^*}$ is a NE of the aggregative game, there exists $\boldsymbol{\chi^*}\in\it{R}^{N(n+1)}$ such that $(\boldsymbol{\chi^*}, \boldsymbol{0}_{Nq})$ is an equilibrium point of \eqref{sf1}.
\end{lemma}

\textbf{Proof.}
{\it Sufficiency:} If $(\boldsymbol{\chi^*},\boldsymbol{\rho^*})$ is the equilibrium point of the system \eqref{sf1}, which satisfies
\begin{subequations}
\begin{align}
\begin{split}
\boldsymbol{0}_{N(n+1)}&=\mathcal{A}\boldsymbol{\chi^*}+\mathcal{B}\boldsymbol{e^*}+\mathcal{D}\boldsymbol{\rho^*}, \label{ep1}
\end{split}\\
\boldsymbol{0}_{Nq}&=(S-K_o\boldsymbol{B}^T\boldsymbol{B}\boldsymbol{U})\boldsymbol{\rho^{*}}.\label{ep2}
\end{align}
\end{subequations}

Denote $\boldsymbol{x}=\col(x_1,\ldots,x_N)$ and $\boldsymbol{\gamma}=[\gamma_1,\ldots,\gamma_N]^T$. By simple derivation, \eqref{ep1} yields that
\begin{subequations}
\begin{align}
\begin{split}
\boldsymbol{0}_{Nn}&=\boldsymbol{H}\boldsymbol{x^*}+\boldsymbol{B}(K_P\boldsymbol{e^*}+K_I\boldsymbol{\gamma^{*}})+\boldsymbol{B}\boldsymbol{U}\boldsymbol{\rho^{*}}, \label{ep21}
\end{split}\\
\boldsymbol{0}_{N\ }&=\boldsymbol{e^*}=-\phi(\boldsymbol{y^*}), \label{ep22}
\end{align}
\end{subequations}
where $\boldsymbol{H}=\blk\{H_1,\cdots,H_N\}$, $K_P=\diag\{k_{P_1},\ldots,k_{P_N}\}$ and $K_I=\diag\{k_{I_1},\ldots,$ $k_{I_N}\}$.
Since $(\boldsymbol{S}-K_o\boldsymbol{B}^T\boldsymbol{B}\boldsymbol{U})$ is Hurwitz, it is derived from \eqref{ep2} that $\boldsymbol{\rho^{*}}=\boldsymbol{0}_{Nq}$. Thus, it is derived from \eqref{ep21} that $\boldsymbol{H}\boldsymbol{x^*}+\boldsymbol{B}K_I\boldsymbol{\gamma^{*}}=\boldsymbol{0}_{Nn}$. From the system \eqref{sf1}, it follows that $\boldsymbol{y^*}=\mathcal{C}\boldsymbol{\chi^*}$. According to \eqref{ep22}, $\phi(\boldsymbol{y^*})=\boldsymbol{0}_{N}$, which satisfies the condition \eqref{psegc}. So, we have the conclusion that $\boldsymbol{y^*}$ is the Nash equilibrium of the game $G=(\mathcal{I},\Omega,\tilde{J})$ .

{\it Necessary:} If $\boldsymbol{y^*}$ is a NE of aggregative game $G=(\mathcal{I},\Omega,\tilde{J})$, we have
\begin{equation*}
\nabla_i{J_i(y_i^*,y_{-i}^*)}=0,   \  \forall{i}\in\mathcal{I},
\end{equation*}
and $\phi(\boldsymbol{y}^*)=\mathbf{0}_N$ such that \eqref{ep22} is satisfied. Then, there exists an  $x_i^*\in\it{R}^{n}$ such that $y_i^*=Cx_i^*$ for all $i\in \mathcal{I}$. Since $(\boldsymbol{S}-K_o\boldsymbol{B}^T\boldsymbol{B}\boldsymbol{U})$ is Hurwitz, $\boldsymbol{\rho^{*}}=\boldsymbol{0}_{Nq}$ and \eqref{ep2} is satisfied. According to that $\gamma_i$ is the integral of $e_i$ for all $i\in \mathcal{I}$ and $\boldsymbol{H}$ is Hurwitz, there exists $\boldsymbol{\gamma^*}\in\it{R}^N$ to satisfy \eqref{ep21}. By the definition of $\boldsymbol{\chi^*}$ which consists of $\boldsymbol{x^*}$ and $\boldsymbol{\gamma^*}$, there exists $\boldsymbol{\chi^*}\in \it{R}^{N(n+1)}$ to satisfy \eqref{ep1}. Therefore, $(\boldsymbol{\chi^*}, \boldsymbol{0}_{Nq})$ is an equilibrium point of
\eqref{sf1}. \qed

Lemma \ref{lemma3} indicates that the NE of aggregative game $G=(\mathcal{I},\Omega,\tilde{J})$ can be reached by the strategy-updating rule \eqref{pi} if the state of system \eqref{sf1} can converge to the equilibrium point. Then, we will analyze the convergence of \eqref{sf1} to illustrate the effectiveness of \eqref{pi}.

\begin{theorem} \label{thm1}
Under Assumptions \ref{as2} and \ref{as3}, if agents with dynamics \eqref{ad1} follow the strategy-updating rule \eqref{pi}, the strategies of all agents can globally asymptotically converge to the unique Nash equilibrium of the game $G=(\mathcal{I},\Omega,\tilde{J})$ with perfect information, provided the following conditions hold

1) $k_{P_i}>0$ and $k_{I_i}>0$ such that $(\mathcal{A}_i,\mathcal{B}_i)$ is controllable and $(\mathcal{A}_i,\mathcal{C}_i)$ is observable, $i\in \mathcal{I}$;

2) disturbance observer gain matrix $k_{o_i}\in \it{R}^p$ such that $S_i-k_{o_i}B^TBU_i$ is Hurwitz, $i\in \mathcal{I}$;

3) the symmetric positive definite matrix $P_i$ is chosen such that $P_i\mathcal{B}_i=\mathcal{C}_i^T$, $i\in \mathcal{I}$.
\end{theorem}

\textbf{Proof.}
Denote $\boldsymbol{\tilde{\chi}}=\boldsymbol{\chi}-\boldsymbol{\chi{^*}}$,  $\boldsymbol{\tilde{y}}=\boldsymbol{y}-\boldsymbol{y^*}$, and $ \boldsymbol{\tilde{e}}=\boldsymbol{e}-\boldsymbol{e^*}$.
Then, the close-loop system is given by
\begin{equation} \label{nc1}
\begin{aligned}
\boldsymbol{\dot{\tilde{\chi}}}&=\mathcal{A}\boldsymbol{\tilde{\chi}}+\mathcal{B}\boldsymbol{\tilde{e}}+\mathcal{D}\boldsymbol{\rho},\\
\boldsymbol{\tilde{y}}&=\mathcal{C}\boldsymbol{\tilde{\chi}},\\
\boldsymbol{\dot{\rho}}&=(\boldsymbol{S}-K_o\boldsymbol{B}^T\boldsymbol{B}\boldsymbol{U})\boldsymbol{\rho}.
\end{aligned}
\end{equation}

For system \eqref{red}, because $H_i$ is Hurwitz, the augmented matrix $\mathcal{A}_i$ has a simple zero eigenvalue and all others with negative real parts. $\rho_i$ is the reference input, which does not affect the poles distribution of the system \eqref{red}. Then, appropriate parameters $k_{P_i}$ and $k_{I_i}$ are chosen to make $(\mathcal{A}_i,\mathcal{B}_i,\mathcal{C}_i)$ be minimal realization. Thus, all poles of the close-loop system \eqref{nc1} have non-positive real parts. From Definition 6.4 in \cite{Khalil.2002}, \eqref{nc1} is positive real. When the observation error is zero or the system is absence of disturbances, i.e., $\boldsymbol{\rho}=\boldsymbol{0}_{Nn}$, the system \eqref{nc1} can be reduced to
\begin{equation} \label{nc12}
\begin{aligned}
\boldsymbol{\dot{\tilde{\chi}}}&=\mathcal{A}\boldsymbol{\tilde{\chi}}+\mathcal{B}\boldsymbol{\tilde{e}},\\
\boldsymbol{\tilde{y}}&=\mathcal{C}\boldsymbol{\tilde{\chi}}.
\end{aligned}
\end{equation}

For the reduced system \eqref{nc12}, $\boldsymbol{\tilde{\chi}}$, $\boldsymbol{\tilde{e}}$, and $\boldsymbol{\tilde{y}}$ are the state, input and output, respectively. All poles of the transfer function matrix $G(s)=\mathcal{C}(sI-\mathcal{A})^{-1}\mathcal{B}$ have non-positive real parts, which can be easily derived from \eqref{red}. Thus, by Lemma \ref{lem6.4},
system \eqref{nc12} is passive. There exists a storage function $V_1=\frac{1}{2}\boldsymbol{\tilde{\chi}}^T\boldsymbol{P}\boldsymbol{\tilde{\chi}}$, where $\boldsymbol{P}=\blk\{P_1,\cdots,P_N\}$, and $P_i=P_i^T>0, \forall i\in\mathcal{I}$. The derivation of $V_1$ along the solutions of \eqref{nc12} is
\begin{equation*}
\begin{split}
\dot{V}_1&= \boldsymbol{\tilde{\chi}}^T\boldsymbol{P}(\mathcal{A}\boldsymbol{\tilde{\chi}}+\mathcal{B}\boldsymbol{\tilde{e}})\\
&=\frac{1}{2}\boldsymbol{\tilde{\chi}}^T(\boldsymbol{P}\mathcal{A}+\mathcal{A}^T\boldsymbol{P})\boldsymbol{\tilde{\chi}}+\boldsymbol{\tilde{\chi}}^T\boldsymbol{P}\mathcal{B}\boldsymbol{\tilde{e}}\\
&=-\frac{1}{2}\boldsymbol{\tilde{\chi}}^T\boldsymbol{Q}^T\boldsymbol{Q}\boldsymbol{\tilde{\chi}}^T+\boldsymbol{\tilde{\chi}}^T\mathcal{C}^T\boldsymbol{\tilde{e}}\\
&\leq
\boldsymbol{\tilde{e}}^T\boldsymbol{\tilde{y}}\\
&=(\boldsymbol{e}-\boldsymbol{e}^*)^T\boldsymbol{\tilde{y}}\\
&=-\big(\phi(\boldsymbol{y})-\phi(\boldsymbol{y^*})\big)^T(\boldsymbol{y}-\boldsymbol{y^*})\\
&\leq-\mu\|\boldsymbol{\tilde{y}}\|^2\leq-\mu \|\mathcal{C}\| \|\boldsymbol{\tilde{\chi}}\|^2.
\end{split}
\end{equation*}
where $\boldsymbol{Q}=\blk\{Q_1,\cdots,Q_N\}$ satisfying $\boldsymbol{P}\mathcal{A}+\mathcal{A}^T\boldsymbol{P}=-\boldsymbol{Q}^T\boldsymbol{Q}$ and $\boldsymbol{P}\mathcal{B}=\mathcal{C}^T$ by Lemma \ref{lem6.2}.
The second inequality derives from the strong monotonicity of function $\phi$ in Assumption \ref{as3}.
It indicates that the system \eqref{nc12} can converge globally exponentially to the origin.

Next, we consider the cascade system \eqref{nc1}. Then, the derivative of $V_1$ along the solution of \eqref{nc1} is
\begin{equation*}
\begin{split}
\dot{V}_1&=\boldsymbol{\tilde{\chi}}^T\boldsymbol{P}(\mathcal{A}\boldsymbol{\tilde{\chi}}+\mathcal{B}\boldsymbol{\tilde{e}}+\mathcal{D}\boldsymbol{\rho})\\
&\leq \boldsymbol{\tilde{e}}^T\boldsymbol{\tilde{y}}+\boldsymbol{\tilde{\chi}}^T\boldsymbol{P}\mathcal{D}\boldsymbol{\rho}\\
&\leq -\mu\|\mathcal{C}\| \|\boldsymbol{\tilde{\chi}}\|^2+\|\boldsymbol{\tilde{\chi}}\| \|\boldsymbol{P}\mathcal{D}\| \|\boldsymbol{\rho}\|,\\
&\leq -\mu(1-b)\|\mathcal{C}\| \|\boldsymbol{\tilde{\chi}}\|^2-\mu b\|\mathcal{C}\| \|\boldsymbol{\tilde{\chi}}\|^2+\|\boldsymbol{\tilde{\chi}}\| \|\boldsymbol{P}\mathcal{D}\| \|\boldsymbol{\rho}\|,\\
\end{split}
\end{equation*}
where $0<b<1$. Then, for $\forall \|\boldsymbol{\tilde{\chi}}\|\geq \frac{\|\boldsymbol{P}\mathcal{D}\|}{\mu b\|\mathcal{C}\|}\|\boldsymbol{\rho}\|$,
\begin{equation*}
\begin{split}
\dot{V}_1&\leq -\mu(1-b)\|\mathcal{C}\| \|\boldsymbol{\tilde{\chi}}\|^2.
\end{split}
\end{equation*}

By Theorem \ref{thm3}, the $\boldsymbol{\tilde{\chi}}$-subsystem in \eqref{nc1} is input-to-state stable. Since $(\boldsymbol{S}-K_o\boldsymbol{B}^T\boldsymbol{B}\boldsymbol{U})$ is Hurwitz, $\boldsymbol{\rho}$ globally asymptotically converges to the origin.  Then, the origin of the cascade system \eqref{nc1} is globally asymptotically stable by Lemma \ref{lemma1}. Therefore, $\boldsymbol{\chi}$ and $\boldsymbol{y}$ of the system \eqref{sf1} globally asymptotically converge to $\boldsymbol{\chi^*}$ and $\boldsymbol{y^*}$, respectively, that is, all players with strategy-updating rule \eqref{pi} asymptotically reach the NE of the game. \qed

\subsection{Strategy-updating Rule with Imperfect Information}
In this section, we consider that the players cannot obtain all players' actions in a large-scale network. So, they need to estimate the aggregator $\sigma$ through communication with their neighbors over a graph $\mathcal{G}$. Denote the estimation of the aggregator $\sigma$ of player $i$ by $\eta_i\in \it{R}$, $\forall i\in \mathcal{I}$. Based on the strategy-updating rule \eqref{pi} and observer \eqref{do}, we employ a dynamic consensus protocol such that $\eta_i=\eta_j$ for all $i,j\in \mathcal{I}$, before all players reach the NE.

Let $e'_i=-F_i(y_i,\eta_i)$ denote the gradient of the player $i$'s cost function coupled with the aggregator estimation $\eta_i$, which is different from the case of perfect information. The strategy-updating rule for the player $i\in \mathcal{I}$ is designed as
\begin{equation}\label{sur2}
\begin{aligned}
\dot{x}_i&=H_ix_i+B_i(k_{P_i}e'_i+k_{I_i}\gamma_i)+B_iU_i\rho_i,\\
y_i &=C_ix_i,
\end{aligned}
\end{equation}
where the disturbance observation error $\rho_i$ is given by
\begin{equation} \label{doe}
\dot{\rho}_i=(S_i-k_{o_i}B_i^TB_iU_i)\rho_i.
\end{equation}

Based on the dynamic consensus protocol in \cite{Ye.2017,Menon.2014}, the estimation of aggregator $\sigma$ is expressed as
\begin{equation} \label{dcp}
\begin{aligned}
\delta\dot{\eta}_i&=-\eta_i-\sum_{j\in\mathcal{N}_i}(\eta_i-\eta_j)-\sum_{j\in\mathcal{N}_i}(\omega_i-\omega_j)+Ny_i, \\
\delta\dot{\omega}_i&=\sum_{j\in\mathcal{N}_i}(\eta_i-\eta_j),
\end{aligned}
\end{equation}
where $\omega_i$ is an auxiliary variable, and $\delta$ is a small positive constant. Because the aggregator estimation is the sum of outputs of all players, it is time-varying. The dynamic average consensus protocol is more suitable than the static one. In order to track the changing aggregator,  $\delta$ is needed to make the protocol \eqref{dcp} a fast subsystem tracking the aggregator in the slow subsystem.

Let $\boldsymbol{e'}=[e'_1,\ldots,e'_N]^T$, $\eta=[\eta_1,\ldots,\eta_N]^T$, and $\omega=[\omega_1,\ldots,\omega_N]^T$. Based on the analysis in the previous section, \eqref{sur2}, \eqref{doe} and \eqref{dcp} can be described as
\begin{equation} \label{sf2}
\begin{aligned}
\boldsymbol{\dot{\chi}}&=\mathcal{A}\boldsymbol{\chi}+\mathcal{B}\boldsymbol{e'}+\mathcal{D}\boldsymbol{\rho},\\
\boldsymbol{y}&=\mathcal{C}\boldsymbol{\chi},\\
\boldsymbol{\dot{\rho}}&=(\boldsymbol{S}-K_o\boldsymbol{B}^T\boldsymbol{B}\boldsymbol{U})\boldsymbol{\rho},\\
\begin{bmatrix}
\rm{\delta\dot{\eta}}\\
\rm{\delta\dot{\omega}}
\end{bmatrix}&=
\begin{bmatrix}
-I-L & -L\\
L & \boldsymbol{0}
\end{bmatrix}
\begin{bmatrix}
\rm{\eta}\\
\rm{\omega}
\end{bmatrix}+
\begin{bmatrix}
N\boldsymbol{y}\\
\boldsymbol{0}_N
\end{bmatrix}.
\end{aligned}
\end{equation}
In the close-loop system \eqref{sf2}, $(\eta,\omega)$ is a fast subsystem, and $(\boldsymbol{\chi,\rho})$ is a slow subsystem. For a connected and undirected graph $\mathcal{G}$ and by using Theorem 5 in \cite{Freeman.2006}, we find that $\eta(t)$ first converges exponentially to equilibrium point $\sum_{i=1}^Ny_i$$\boldsymbol{1_N}$ for slowly changing $\boldsymbol{y}$. Then, $(\boldsymbol{\chi},\boldsymbol{\rho})$ asymptotically converge to its equilibrium point $(\boldsymbol{\chi^*},\boldsymbol{0}_{Nq})$.

\begin{theorem} \label{thm2}
Under Assumptions \ref{as2}-\ref{as3} and the conditions in Theorem \ref{thm1}, for a connected and undirected graph $\mathcal{G}$, there exists a positive constant $\delta^*$ such that for every $0<\delta<\delta^*$, $\eta(t)$ exponentially converges to $\sum_{i=1}^Ny_i\boldsymbol{1_N}$. Agents with dynamics \eqref{ad1} follow the strategy-updating rule \eqref{sur2}. Then, the strategies of all agents can globally asymptotically converges to the unique Nash equilibrium $\boldsymbol{y^*}$ of aggregative game $G=(\mathcal{I},\Omega,\tilde{J})$ with imperfect information.
\end{theorem}

\textbf{Proof.}
Denote $\boldsymbol{\tilde{\chi}}=\boldsymbol{\chi}-\boldsymbol{\chi^*}$,  $\boldsymbol{\tilde{y}}=\boldsymbol{y}-\boldsymbol{y^*}$, $\boldsymbol{\tilde{e}'}=\boldsymbol{e'}-\boldsymbol{e'^*}$, $\tilde{\eta}=\eta-\bar{\eta}$, and $\tilde{\omega}=\omega-\bar{\omega}$. Then, system \eqref{sf2} can be rewritten as
\begin{equation} \label{nc2}
\begin{aligned}
\boldsymbol{\dot{\tilde{\chi}}}&=\mathcal{A}\boldsymbol{\tilde{\chi}}+\mathcal{B}\boldsymbol{\tilde{e}'}+\mathcal{D}\boldsymbol{\rho},\\
\boldsymbol{\tilde{y}}&=\mathcal{C}\boldsymbol{\tilde{\chi}},\\
\boldsymbol{\dot{\rho}}&=(\boldsymbol{S}-K_o\boldsymbol{B}^T\boldsymbol{B})\boldsymbol{U})\boldsymbol{\rho},\\
\begin{bmatrix}
\rm{\delta\dot{\tilde{\eta}}}\\
\rm{\delta\dot{\tilde{\omega}}}
\end{bmatrix}&=
\begin{bmatrix}
-I-L & -L\\
L & \boldsymbol{0}
\end{bmatrix}
\begin{bmatrix}
\rm{\tilde{\eta}}\\
\rm{\tilde{\omega}}
\end{bmatrix},
\end{aligned}
\end{equation}
where $\bar{\eta}$ and $\bar{\omega}$ are the quasi-steady states of $\eta$ and $\omega$ for fixed $\boldsymbol{y}$, respectively. According to the stability analysis of singularly perturbed system, there are three steps as follows.

1) Quasi-steady state analysis: Let $\delta=0$ freeze $\eta$ and $\omega$ at the quasi-steady states, i.e., $\eta_i=\bar{\eta}_i=\sum_{j=1}^N y_i$ and $\omega_i=\bar{\omega}_i,  \forall{i\in \mathcal{I}}$. Then $\boldsymbol{\tilde{e}'}=-F(\boldsymbol{\tilde{y}}+\boldsymbol{y^*},\sigma)=-\phi(\boldsymbol{\tilde{y}}+\boldsymbol{y^*})$, the system \eqref{nc2} can be reduced as \eqref{nc1}. According to the analysis in Theorem \ref{thm1}, the reduced system can be globally asymptotically stable at $(\boldsymbol{\chi^*,0}_{Nq})$.

2) Boundary-layer analysis: Let $\tau=t/\delta$, and $\delta=0$ freeze $\boldsymbol{y}$. The boundary-layer system in \eqref{nc2} is rewritten in $\tau$-time scale as
\begin{equation} \label{bl}
\begin{bmatrix}
\rm{\frac{d\tilde{\eta}}{d\tau}}\\
\rm{\frac{d\tilde{\omega}}{d\tau}}
\end{bmatrix}=
\begin{bmatrix}
-I-L & -L\\
L & \boldsymbol{0}
\end{bmatrix}
\begin{bmatrix}
\rm{\tilde{\eta}}\\
\rm{\tilde{\omega}}
\end{bmatrix}.
\end{equation}

The following analysis method is similar to \cite{Ye.2017,Menon.2014}. Let $r$ be an $N$ dimensional column vector such that $r^TL=0$. $R=[\tilde{R}, r]$ denotes an $N\times{N}$ dimensional orthonormal matrix. Then, let $\tilde{\omega}=R\big[\begin{smallmatrix}\tilde{\omega}'\\ \tilde{\omega}_0\end{smallmatrix}\big]$, which is decomposed into an $N-1$ dimensional column vector $\tilde{\omega}'$ and a scalar $\tilde{\omega}_0$. Because $\tilde{\omega}_0$ does not interact with other states in the boundary-layer system, \eqref{bl} can be written as
\begin{equation*}
\begin{bmatrix}
\rm{\frac{d\tilde{\eta}}{d\tau}}\\
\rm{\frac{d\tilde{\omega}'}{d\tau}}
\end{bmatrix}=
\begin{bmatrix}
-I-L & -L\tilde{R}\\
\tilde{R}^TL & \boldsymbol{0}
\end{bmatrix}
\begin{bmatrix}
\rm{\tilde{\eta}}\\
\rm{\tilde{\omega}'}
\end{bmatrix}.
\end{equation*}
For a connected and undirected graph $\mathcal{G}$ and from Lemma 2.2 in \cite{Menon.2014}, $\big[\begin{smallmatrix}
-I-L & -L\tilde{R}\\
\tilde{R}^TL & \boldsymbol{0}
\end{smallmatrix}\big]$ is Hurwitz. $(\tilde{\eta},\tilde{\omega}')$ can exponentially converge to the origin, uniformly in $(t,\boldsymbol{y})$. Based on the converse Lyapunov theorem in \cite[Theorem 4.14]{Khalil.2002}, there exists a Lyapunov function $V_2$ for the boundary-layer system \eqref{bl}. It satisfies $c_1(\|\tilde{\eta}\|^2+\|\tilde{\omega}\|^2)\leq{V_2}\leq c_2(\|\tilde{\eta}\|^2+\|\tilde{\omega}\|^2)$ and $\dot{V_2}\leq -c_3(\|\tilde{\eta}\|^2+\|\tilde{\omega}\|^2)$ for some positive constants $c_1,c_2$ and $c_3$.

3) Comprehensive analysis: Since $\boldsymbol{y}$ is not a true constant parameter, we have to keep track of the effect of the interconnection between the slow and fast dynamics. Now consider the composite Lyapunov function candidate
\begin{equation*}
V=(1-\varepsilon)V_1+\varepsilon{V_2}, \  0<\varepsilon<1,
\end{equation*}
where $\varepsilon$ is a constant. The derivative of $V$ along the trajectory of the full system \eqref{nc2} is
\allowdisplaybreaks[4]
\begin{align*}
\dot{V}&=(1-\varepsilon)\dot{V}_1+\varepsilon(\frac{1}{\delta}-1)\dot{V}_2+(1-\varepsilon)\frac{\partial{V}_1}{\partial{\boldsymbol{\tilde{\chi}}}}(\boldsymbol{\dot{\tilde{\chi}}}\mid_{\tilde{\eta}}
       -\boldsymbol{\dot{\tilde{\chi}}}\mid_{\tilde{\eta}=0})\\
       &\leq -(1-\varepsilon)\mu(1-b)\|\mathcal{C}\|\|{\boldsymbol{\tilde{\chi}}}\|^2-\varepsilon(\frac{1}{\delta}-1)c_3(\|{\tilde{\eta}}\|^2+\|\tilde{\omega}\|^2)\\
       &\ \ \ +(1-\varepsilon)\boldsymbol{\tilde{\chi}}^T\boldsymbol{P}\mathcal{B}(\boldsymbol{\tilde{e}'}\mid_{\tilde{\eta}}-\boldsymbol{\tilde{e}'}\mid_{\tilde{\eta}=0})\\ &= -(1-\varepsilon)\mu(1-b)\|\mathcal{C}\|\|{\boldsymbol{\tilde{\chi}}}\|^2 -\varepsilon(\frac{1}{\delta}-1)c_3(\|{\tilde{\eta}}\|^2+\|\tilde{\omega}\|^2)\\
       &\ \ \ +(1-\varepsilon)\boldsymbol{\tilde{\chi}}^T\boldsymbol{P}\mathcal{B} (-F(\boldsymbol{\tilde{y}+y^*},\tilde{\eta}+\bar{\eta})+F(\boldsymbol{\tilde{y}+y^*},\bar{\eta}))\\
       &\leq -(1-\varepsilon)\mu(1-b)\|\mathcal{C}\|\|{\boldsymbol{\tilde{\chi}}}\|^2 -\varepsilon(\frac{1}{\delta}-1)c_3(\|{\tilde{\eta}}\|^2+\|\tilde{\omega}\|^2)\\
       &\ \ \ +(1\!-\!\varepsilon)\|\boldsymbol{\tilde{\chi}}\|\|\boldsymbol{P}\mathcal{B}\| \|{F(\boldsymbol{\tilde{y}\!+\!y^*}\!,\!\tilde{\eta}\!+\!\bar{\eta})\!-\!F(\boldsymbol{\tilde{y}+y^*},\bar{\eta})}\|\\
       &\leq -(1-\varepsilon)\mu(1-b)\|\mathcal{C}\|\|{\boldsymbol{\tilde{\chi}}}\|^2 -\varepsilon(\frac{1}{\delta}-1)c_3(\|{\tilde{\eta}}\|^2+\|\tilde{\omega}\|^2)\\ &\ \ \ +\theta(1-\varepsilon)\|\boldsymbol{P}\mathcal{B}\|\|\boldsymbol{\tilde{\chi}}\|\|\tilde{\eta}\| \\
       &=-(1-\varepsilon)\begin{bmatrix}
\|{\boldsymbol{\tilde{\chi}}}\|\\ \|{\tilde{\eta}}\|
\end{bmatrix}^T
\Lambda \begin{bmatrix}
\|{\boldsymbol{\tilde{\chi}}}\|\\ \|{\tilde{\eta}}\|
\end{bmatrix}-\varepsilon(\frac{1}{\delta}-1)c_3\|\tilde{\omega}\|^2,
\end{align*}
where
\begin{equation*}
\Lambda=
\begin{bmatrix}
\mu(1-b)\|\mathcal{C}\| & -\frac{1}{2}\theta\|\boldsymbol{P}\mathcal{B}\|\\
-\frac{1}{2}\theta\|\boldsymbol{P}\mathcal{B}\| & \frac{\varepsilon(1-\delta)}{(1-\varepsilon)\delta}c_3
\end{bmatrix}.
\end{equation*}
 The sufficient condition for $\dot{V}<0$ is that $\Lambda$ is positive definite matrix, i.e., there exists a positive constant $\delta^*$, for all $0<\delta<\delta^*$, $(\boldsymbol{\tilde{\chi}},\tilde{\eta})$ asymptotically converges to the origin of \eqref{nc2} \cite[Theorem 11.3]{Khalil.2002}, which indicates system \eqref{sf2} can be stable asymptotically at the Nash equilibrium $\boldsymbol{y*}$.  \qed

\begin{remark}
In recent work \cite{Zhang.2019}, aggregative games of multiple complex dynamic systems with external disturbances were studied. However, the devised algorithm was hard to be applied directly here due to the following three aspects. First, the dynamics of players considered in \cite{Zhang.2019} are nonlinear and with unity relative degrees while in this paper general linear systems are considered without the assumption of relative degrees. Second, the external disturbances considered in \cite{Zhang.2019} are generated by unknown systems, while a class of deterministic disturbances are considered in this paper which can be compensated by observers. Third, the dynamic average consensus protocol employed in \cite{Zhang.2019} requires that initial states be located at the origin, while the consensus protocol proposed in this paper does not require any initial conditions.
\end{remark}

\begin{remark}
By the idea of the singular perturbation method in \cite[Chapeter 11]{Khalil.2002}, the strategy-updating rule designed in this section can be formulated as a singular perturbation model. To make \eqref{dcp} a fast system, the parameter $\delta$ is designed in \eqref{dcp}.  It is convenient to analyze the performance of the network of general linear systems by designing a small parameter $\delta$. Furthermore, the fast system \eqref{dcp} is easily realized by the embedded technology. The parameter may be obtained in a distributed way, which is a problem to be solved in the future.
\end{remark}

\section{Simulations}
In this section, we will present two examples to verify the effectiveness of the proposed strategy-updating rules for all players in aggregative games.
\subsection{Network of Double-integrator Agents}
\begin{figure}
\begin{center}
\includegraphics[width=8cm]{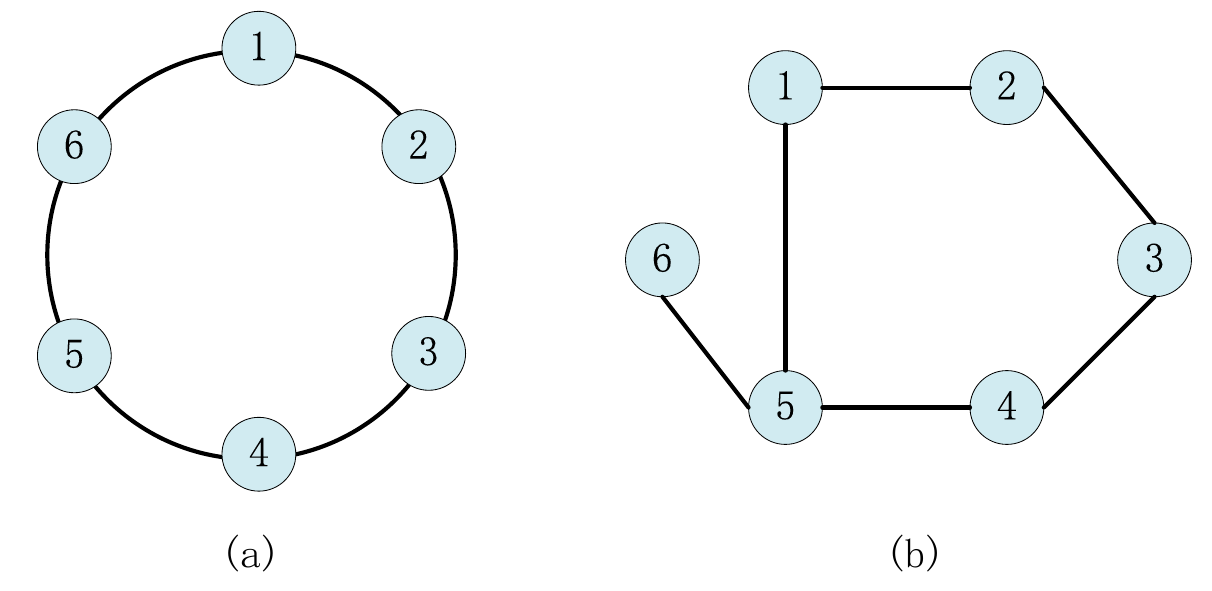}
\caption{Communication graphs.}
\label{fig.4}
\end{center}
\end{figure}

Consider a multi-agent system with six double-integrator agents described by \eqref{ad1}, whose communication topology is depicted in Fig. \ref{fig.4} (a), where $A=\big[\begin{smallmatrix}0 & 1\\0 & 0\end{smallmatrix}]$, $B=[\begin{smallmatrix} 0\\1\end{smallmatrix}]$, and $C=[\begin{smallmatrix}1 & 0\end{smallmatrix}]$. Agent $i$'s cost function is given by
\begin{equation} \label{cfe1}
J_i(y_i,y_{-i})=a_i(y_i-b_i)^2-(c_0-C\sigma)y_i,
\end{equation}
where $a=[a_1,\ldots,a_6]^T=[1.0, 0.5, 0.8, 0.7, 1.1, 0.6]^T$, $b=[b_1,\ldots,b_6]^T=[6, 10, 7, 8, 6, 12]^T$, $c_0=50$, $C=0.5$. According to the conditions stated in Theorem \ref{thm1}, we design $K_i=[1,10]$, $k_{P_i}=12$, $k_{I_i}=2$, and $k_{o_i}=[12,15,12,15,12, 15]^T$, $\forall i\in \mathcal{I}$. Besides, $\delta=0.1$ in the games with incomplete information. Disturbance $d_i$ is a sinusoidal signal which acts on the control channel.  It is clear that cost functions \eqref{cfe1} and \eqref{cfe2} satisfy the conditions in Assumptions \ref{as2}-\ref{as3}.

\begin{figure}
\begin{center}
\includegraphics[width=7cm]{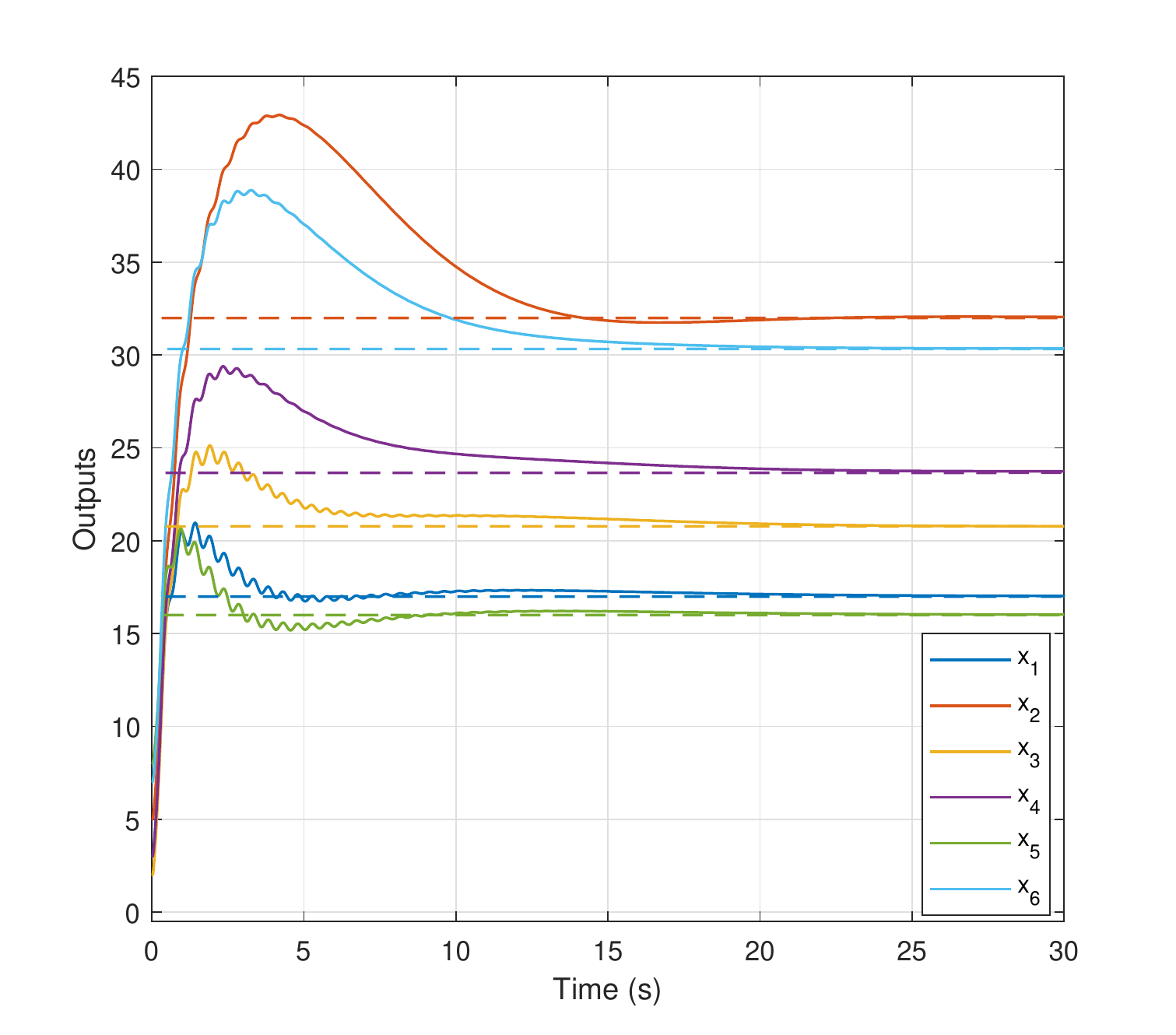}
\caption{Outputs of double-integrator agents by \eqref{pi}.}
\label{fig2}
\end{center}
\end{figure}

\begin{figure}
\begin{center}
\includegraphics[width=7cm]{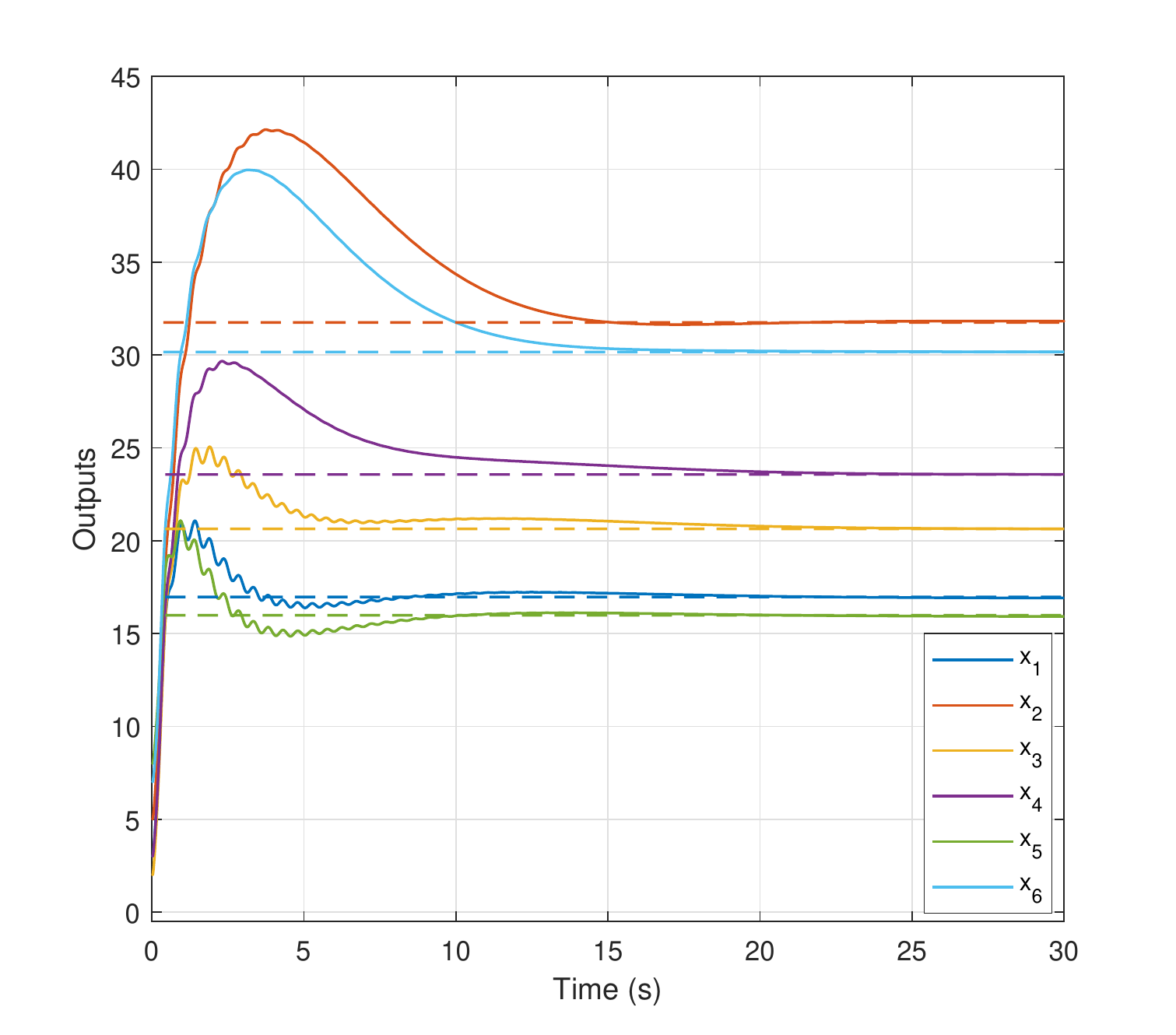}
\caption{Outputs of double-integrator agents by \eqref{sur2}.}
\label{fig3}
\end{center}
\end{figure}

The simulation results of strategy-updating rules \eqref{pi} and \eqref{sur2} for the double-integrator agents with perfect and imperfect information are shown in Figs. \ref{fig2} and \ref{fig3}, respectively. It can be seen that the multi-agent system with agent dynamics \eqref{ad1} under the two strategy-updating rules achieves the NE asymptotically. Moreover, the strategy-updating rule \eqref{sur2} in the case of imperfect information has the same convergence performance as that with perfect information.

\begin{table}[!t]
\renewcommand{\arraystretch}{1.3}
\caption{Description of Variables and Parameters Appearing in the Turbine-generator Model}
\label{Table 1}
\centering
\begin{tabular}{ll}
\hline\hline \\[-3mm]
\bfseries & \bfseries state variables\\
\hline
$P_i$ & output power, in p.u.\\
$X_{ei}$ & valve opening, in p.u.\\
$w_i$ & relative speed, in rad/s\\
\hline
\bfseries  & \bfseries parameters\\
\hline
$T_{mi}$ & time constants of machine's turbine, in s\\
$K_{mi}$ & gain of machine's turbine\\
$T_{ei}$ & time constants of speed governor, in s\\
$K_{ei}$ & gain of speed governor\\
$R_i$ & regulation constant of machine's turbine, in p.u.\\
$w_0$ & synchronous machine speed, in rad/s\\
$D_i$ & per unit damping constant\\
$H_i$ & inertia constant, in s\\
\hline
$u_i$ & control input of generator system, in p.u.\\
\hline\hline
\end{tabular}
\end{table}

\begin{table}
\renewcommand{\arraystretch}{1.2}
\caption{System parameters}
\scalebox{0.68}{
\label{Table 2}
\centering
\begin{tabular}{cccccccccccccc}
\hline\hline \\[-3mm]
~ & $T_{mi}$ & $T_{ei}$ & $K_{mi}$ & $K_{ei}$ & $D_i$ & $H_i$ & $R_i$ & $\alpha_i$ & $\beta_i$ & $\xi_i$ & $P_i(0)$ & $X_{ei}(0)$ & $w_i(0)$\\
\hline
G$\sharp$1 & 0.35 & 0.10 & 1.0 & 1.0 & 5.0 & 4.0 & 0.05 & 5 & 12 & 1.0 & 30 & 6 & 4.3\\
G$\sharp$2 & 0.30 & 0.12 & 1.1 & 1.1 & 4.0 & 3.5 & 0.04 & 8 & 10 & 0.5 & 25 & 5 & 3.5\\
G$\sharp$3 & 0.28 & 0.08 & 0.9 & 0.9 & 3.0 & 2.8 & 0.03 & 6 & 11 & 0.8 & 20 & 4 & 3.0\\
G$\sharp$4 & 0.40 & 0.11 & 1.2 & 1.2 & 4.5 & 4.2 & 0.06 & 9 & 11 & 0.7 & 35 & 7 & 4.8\\
G$\sharp$5 & 0.43 & 0.90 & 0.8 & 0.8 & 3.5 & 3.0 & 0.04 & 7 & 13 & 1.1 & 28 & 5 & 4.0\\
G$\sharp$6 & 0.35 & 0.10 & 1.0 & 1.0 & 5.0 & 4.0 & 0.05 & 8 & 14 & 0.6 & 37 & 8 & 5.0\\
\hline\hline
\end{tabular}}
\end{table}

\subsection{Generation Systems in Electricity Markets}
In the electricity market, power plants as participants compete with each other to obtain generation index which minimize their own costs. It can be described by an aggregative game. We consider a network of six generation systems communicating with each other via an undirected and connected graph as shown in Fig. \ref{fig.4} (b). The cost function of the generation system $i$ is described by
\begin{equation} \label{cfe2}
\begin{split}
J_i(P_i,P_{-i})&=c_i(P_i)-p(\sigma)P_i\\
               &=\alpha_i+\beta_iP_i+\xi_iP_i^2-(p_0-a\sigma)P_i,
\end{split}
\end{equation}
where $P_i\in\it{R}$ is the output power of the generation system $i$, in $p.u.$, $P_{-i}=[P_1,\ldots,P_{i-1},P_{i+1},\ldots,P_N]^T$, $c_i(P_i)$ is the generation cost, $p(\sigma)$ is electricity price, $\alpha_i,\beta_i$, and $\xi_i$ are characteristics of the generation system $i$, $p_0$ and $a$ are constants, and $\sigma=\sum_{i=1}^N P_i$ denotes the aggregator.

Without regard to the mechanical and electromagnetic loss, the classical dynamical model of the $i$th turbine-generator is governed by (refer to \cite{Bergen.1986,Kundur.1994})
\begin{equation*}
\left.\begin{aligned}
\begin{bmatrix}
\dot{P_i} \\ \dot{X_{ei}} \\ \dot{w_i}
\end{bmatrix}
&=\begin{bmatrix}
-\frac{1}{T_{mi}} & \frac{K_{mi}}{T_{mi}} & 0\\
0 & -\frac{1}{T_{ei}} & -\frac{K_{ei}}{T_{ei}R_iw_0}\\
0 & 0 & -\frac{D_i}{2H_i}
\end{bmatrix}
\begin{bmatrix}
P_i \\ X_{ei} \\ w_i
\end{bmatrix}\\
&\ \ \ +\begin{bmatrix}
0 \\ \frac{1}{T_{ei}} \\ 0
\end{bmatrix}
(u_i+d_i),\\
P_i&=\begin{bmatrix}
1 & 0 & 0\\
\end{bmatrix}
\begin{bmatrix}
P_i \\ X_{ei} \\ w_i
\end{bmatrix},
\end{aligned}\right.
\end{equation*}
where the notation for the generator system model is given in Table \ref{Table 1}. The persistent disturbance is modeled by $d_i=m_i\sin{(2\pi ft)}$, where $m_i$ is an unknown amplitude and $f=500Hz$. The parameters of the six generation systems are brought from \cite{Deng.2019,Guo.2000} and given in Table \ref{Table 2}.  The parameters in strategy-updating rules \eqref{pi} and \eqref{sur2} are given by $\delta=0.1$, $k_{P_i}=k_{I_i}=1$, $K_i=[1,1.2,0.8,1.1,0.9,1]$, and $k_{o_i}=[4,4,4,4,4,8;4,4,4,4,4,8]^T$, which satisfy the conditions in Theorems 2 and 3.
\begin{figure}
\begin{center}
\includegraphics[width=7cm]{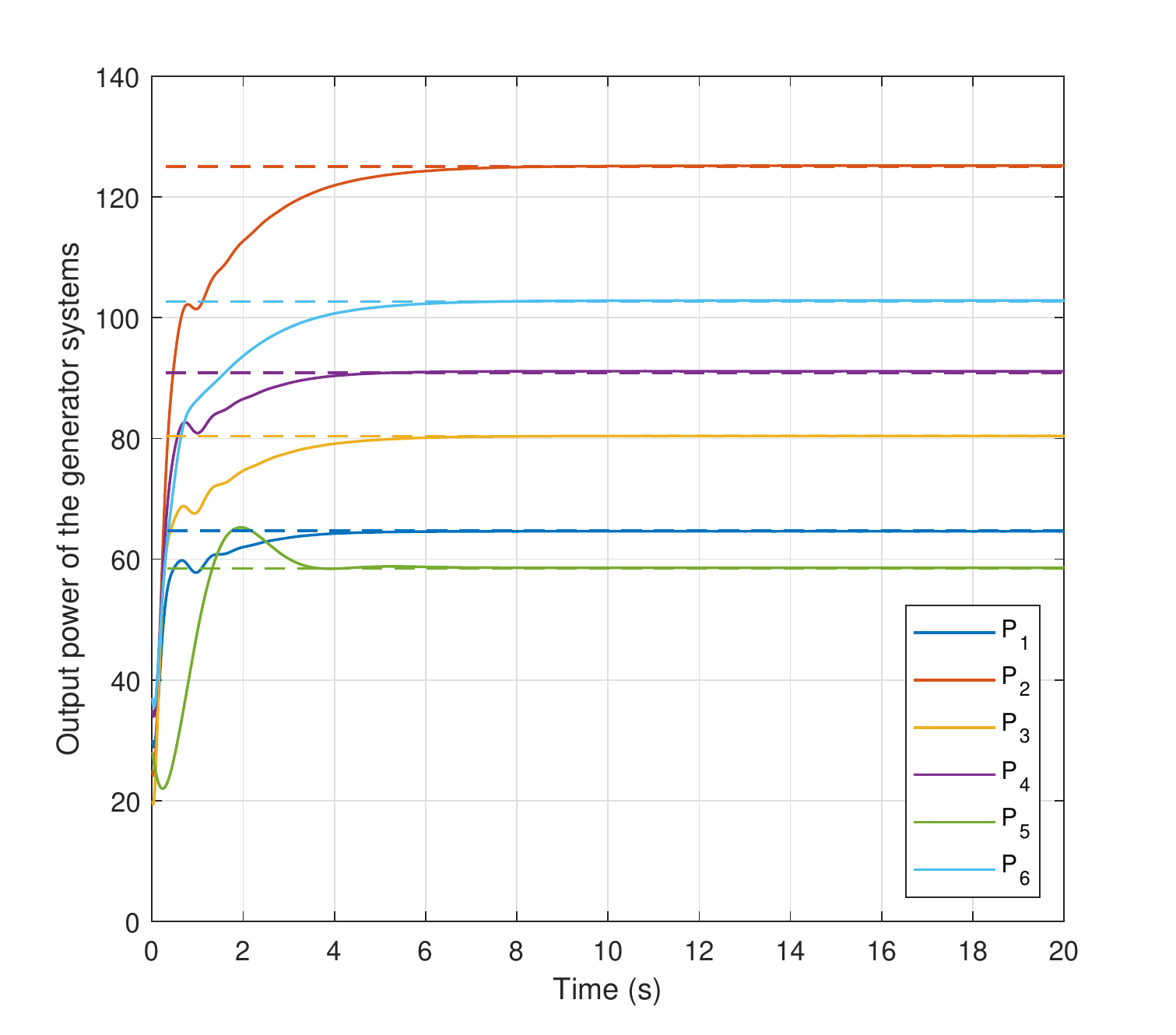}
\caption{Output power of generator systems by strategy-updating rule \eqref{pi}.}
\label{fig.5}
\end{center}
\end{figure}

\begin{figure}
\begin{center}
\includegraphics[width=7cm]{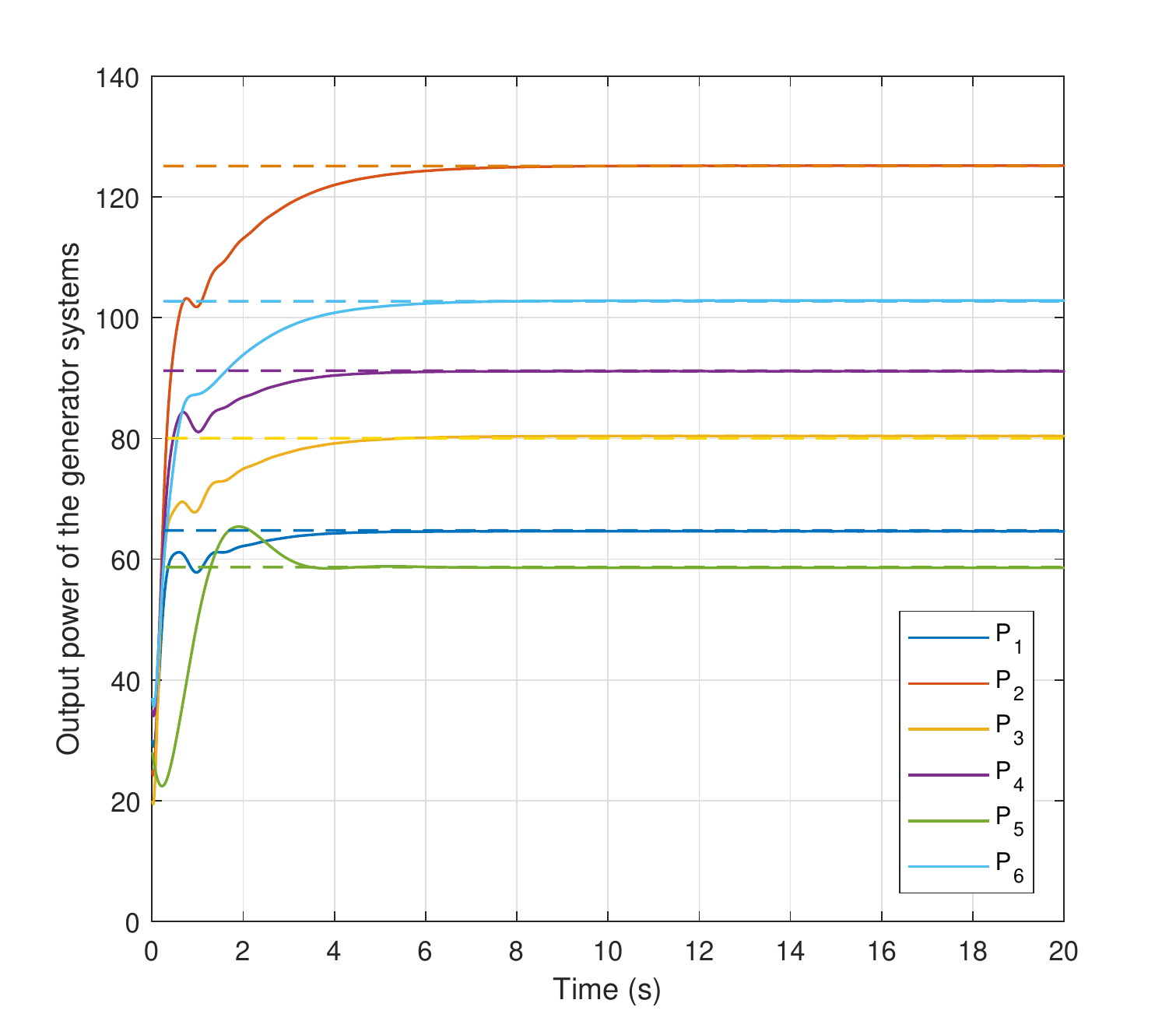}
\caption{Output power of generator systems by strategy-updating rule \eqref{sur2}.}
\label{fig.6}
\end{center}
\end{figure}

\begin{figure}
\begin{center}
\includegraphics[width=7cm]{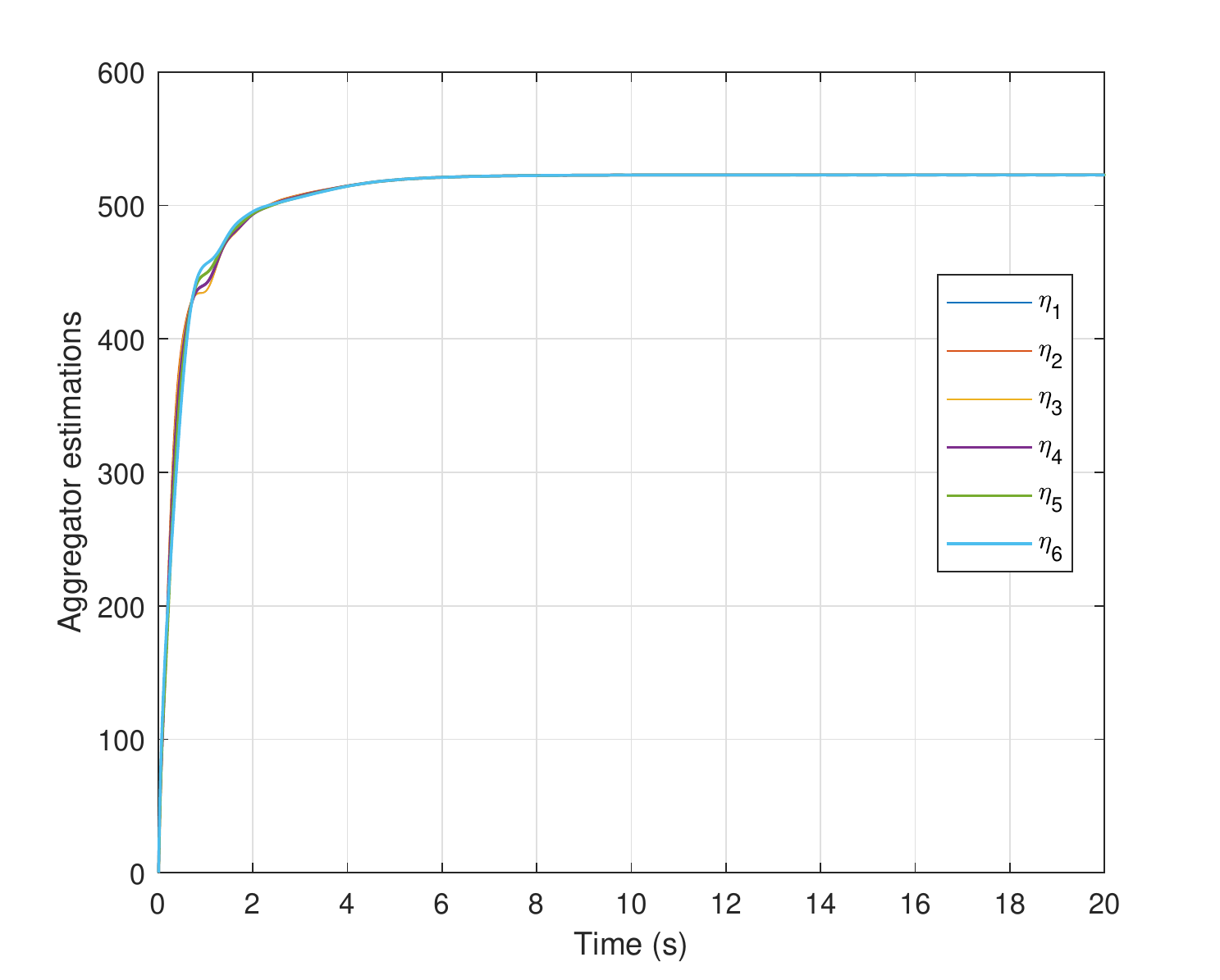}
\caption{Aggregator estimation of six generator systems.}
\label{fig.7}
\end{center}
\end{figure}

In the perfect information case, the regulation process of output powers of the six generators by the strategy-updating rule \eqref{pi} is shown in Fig.~\ref{fig.5}. It is easily seen that the strategies of the six generators can arrive at the Nash equilibrium of the aggregative game. The evolution of strategies of the six generators by updating strategies with the rule \eqref{sur2} in the imperfect information setting is shown in Fig.~\ref{fig.6}. We can see that the strategies of the six generators can reach the Nash equilibrium by only local interactions with neighbors. Fig.~\ref{fig.7} shows that the aggregator estimations of generator systems rapidly reach a consensus. Figs.~\ref{fig.6}-\ref{fig.7} illustrate that the two-time scale designed in the proposed rule \eqref{sur2} can facilitate strategy-updating in the imperfect information setting.

\section{Conclusion}
Aggregative games of multiple general linear systems subject to persistence disturbance were considered in this paper. Based on the internal model, two strategy-updating rules have been proposed for the players with perfect or imperfect information to reject disturbances.  Besides, the integral of cost functions' gradient was introduced in strategy-updating rules to overcome the difficulty that gradient dynamics cannot be directly applied to general linear systems. The convergence of the proposed rules are analyzed via passive theory, singular perturbation analysis and Lyapunov stability theory. All players, who update their strategies according to the two rules, can asymptotically converge to the Nash equilibrium, which was verified by simulation examples.

Some future research directions are considered as follows.
1) The practical dynamic systems, such as hybrid systems \cite{Zheng.2018}, switched systems \cite{Li.2017,Sun.2017}, could be considered in the game.
2) In cyber-physical systems, agents dynamics may be continuous or discrete-time, and information among agents is exchanged in a discrete-time manner. Thus, the mechanism is worth studying in a discrete-time setup.
3) Private constraints of each player and coupled constraints may be considered in decision spaces of the game. Moreover, nonsmooth or non-convex functions could be considered to generalize the results.

\ifCLASSOPTIONcaptionsoff
  \newpage
\fi

\end{document}